\documentclass{amsart}

\usepackage{latexsym,amssymb,amsmath}
\usepackage{fullpage}

\usepackage{eufrak}

\newtheorem{theorem}{Theorem}
\newtheorem{lemma}{Lemma}
\newtheorem{proposition}{Proposition}
\newtheorem{corollary}{Corollary}
\newtheorem{definition}{Definition}
\newtheorem{example}{Example}
\newtheorem{remark}{Remark}

\def\a{\alpha}

\def\c{\chi}
\def\g{\gamma}

\def\l{\lambda}
\def\L{\Lambda}
\def\m{\mu}
\def\n{\nu}

\def\s{\sigma}

\def\x{{\bf x}}

\def\S{\EuFrak{S}}

\usepackage{youngtab}
\newcommand{\nil}{\mbox{}}
\def\ten{10}

\def\twelve{12}

\title{Combinatorial operators for Kronecker powers of representations of $\S_n$}

\author{Alain Goupil$^{1}$ and Cedric Chauve$^{2}$} 
\address{ C\'egep du Vieux Montr\'eal and LaCIM, UQ\`AM,  CP 8888,
  succ. Centre-ville H3C 3P8, Montr\'eal, QC, Canada} 
\email{goupil@math.uqam.ca}
\address{Computer Science Department and LaCIM, UQ\`AM, CP 8888,
  succ. Centre-ville H3C 3P8, Montr\'eal, QC, Canada}
\email{chauve@lacim.uqam.ca}

\begin{document}

\footnotetext[1]{Work partially supported by a grant from NSERC.}
\footnotetext[2]{Work partially supported by a grant from UQ\`AM.}

\begin{abstract}
  We present  combinatorial operators for the expansion of the
  Kronecker product of irreducible representations of the symmetric
  group $\S_n$.  
  These combinatorial operators are defined in the ring of symmetric
  functions and act on the Schur functions basis. 
  This leads to a combinatorial description of the Kronecker
  powers of the irreducible representations indexed with the partition
  $(n-1,1)$ which specializes  the concept of oscillating tableaux in
  Young's lattice previously defined by S. Sundaram.  
  We call our specialization {\it Kronecker tableaux}. 
  Their combinatorial analysis leads to enumerative results for the multiplicity
  of any irreducible representation in the Kronecker powers of the form
  ${\c^{(n-1,1)}}^{\otimes k}$.     
\end{abstract}

\maketitle

\section{Introduction}
\label{sec:intro}

The subject of the present work is the investigation of the Kronecker
product, sometimes called inner tensor product, of irreducible
representations of the symmetric group $\S_n$. Given two linear
representations  
\begin{eqnarray*}
  {\begin{subarray}{c} A:\S_n\rightarrow Aut(V)\\
      \ \s\,\mapsto\, A(\s)
    \end{subarray}}
  \qquad
  {\begin{subarray}{c} B:\S_n\rightarrow Aut(W)\\
      \ \s\,\mapsto\, B(\s)
    \end{subarray}}
\end{eqnarray*}
which associate linear operators to permutations $\s\in\S_n$, the
Kronecker product of $A$ and $B$, denoted $A\otimes B$, is the
representation of $\S_n$ defined by 
\begin{eqnarray*}
  {\begin{subarray}{c} A\otimes B:\S_n\rightarrow Aut(V\otimes W)\\ 
      \phantom{A\otimes B \S} \s\,\mapsto\, A(\s)\otimes B(\s)
    \end{subarray}}
\end{eqnarray*} 
which is the action on the tensor product $V\otimes W$ of vector
spaces  $V$ and $W$ by means of the tensor product  $A(\s)\otimes
B(\s)$ of the linear operators  $A(\s)$ and  $B(\s)$.   

The irreducible representations of $\S_n$ are the representations
which are indecomposable as direct sums of representations.
They are indexed with the partitions $\l$ of $n$:
$$A^\l:\S_n\rightarrow Aut(V), \s\mapsto A^\l(\s).$$  
The Kronecker product $A^\l\otimes A^\m$ of two  irreducible
representations of $\S_n$ is in general not an irreducible
representation of $\S_n$ and the fundamental problem of expanding it
as a direct sum of irreducible representations   
\begin{eqnarray*}
  A^\l\otimes A^\m=\sum_{\a\vdash n}t_{\l,\m}^\a A^\a
\end{eqnarray*} 
goes  back  to the foundation  of the representation theory. This
problem was studied by Murnaghan \cite{Mu,Mu2}, Littlewood \cite{Li}
and more recently by Lascoux \cite{La}, Garsia and Remmel \cite{GR},
Thibon {\em et al.} \cite{STW}, and others (see \cite{BK} and
references therein).   

To obtain the decomposition coefficients $t_{\l,\m}^\a$, one can use
the characters of the corresponding representations.   
The character of a representation $A$ of $\S_n$ is the map $\c^A$
which sends permutations $\s\in \S_n$ to  the traces of $A(\s)$:  
\begin{eqnarray*}
  {\begin{subarray}{c} \c^A:\S_n\rightarrow \mathbb{C}
      \\
      \phantom{A\otimes B BBBB} \s\,\mapsto\, tr(A(\s))
    \end{subarray}
  }
\end{eqnarray*} 
The fact that the character of the Kronecker product $A\otimes B$ of
two representations is obtained by multiplying the traces of the
linear operators $A(\s)$ and $B(\s)$: 
\begin{eqnarray}
  \label{eq1}
  {\begin{subarray}{c} \c^{A\otimes B}:\S_n\rightarrow \mathbb{C}
      \\
      \phantom{A\otimes B \S A\otimes B \S BBBBBB }  \s\,\mapsto\,
      tr(A(\s)) tr( B(\s)) 
    \end{subarray}
  }
\end{eqnarray} 
is a straightforward consequence of the definition of tensor product
of linear operators.  
Hence one can use property (\ref{eq1}) of characters and the character
table, indexed by integer partitions 
of $n$ (see Table \ref{tab:char} for example) to compute the
characters of a Kronecker product of two irreducible representations.  
Let us recall that the character values $\c^A(\s_{1}),\,\c^A(\s_{2})$  on two
  permutations $\s_1,\s_2\in C_\m$ in the same conjugacy class  $C_\m$
  are equal. Therefore we use the notation $\c^\l_\m$ for the value of
  the irreducible character $\c^\l$ on any element of the conjugacy
  class $C_\m$.    

\begin{example} \label{ex1} \em
  Let Table \ref{tab:char} be the table of the irreducible characters
  of $\S_4$. 
                                %
                                %
\begin{table}[htbp]
  \small
  \begin{center}
    \begin{tabular}{c|r|r|r|r|c}
      {$\l\setminus$}{$\m$} &{ $(4)$} & {$(3,1)$}
      & {$(2,2)$} & {$(2,1^2)$} & {$(1^4)$}
      \rule[-.16cm]{0cm}{.5 cm}\\
      \hline {$(4)$}  &{ $1$} & {$1$} & {$1$} & {$1$} & {$1$}
      \rule[-.16cm]{0cm}{.5 cm}\\ \hline
      {$(3,1)$}  &{ $-1$} & {$0$} & {$-1$} & {$1$} & {$3$}
      \rule[-.16cm]{0cm}{.5 cm}\\ \hline
      {$(2,2)$} &{ $0$}  &{ $-1$} &{ $2$} & {$0$} & {$2$}
      \rule[-.16cm]{0cm}{.5 cm}\\ \hline
      {$(2,1^2)$}  &{ $1$} & {$0$} & {$-1$} & {$-1$} & {$3$}
      \rule[-.16cm]{0cm}{.5 cm}\\ \hline
      {$(1^4)$} &{ $-1$} & {$1$} & {$1$} & {$-1$} & {$1$}
      \rule[-.16cm]{0cm}{.5 cm}
    \end{tabular}
    \smallskip
    \caption{Irreducible characters $\chi^\l_\m$ of $S_4$. \label{tab:char}}
  \end{center}
\end{table}
 
 The character of the Kronecker product $A^{(3,1)}\otimes
  A^{(3,1)}$ of the irreducible representation  $A^{(3,1)}$ with
  itself  is obtained by multiplying each element of the row-vector 
  $(-1,0,-1,1,3)$ in Table \ref{tab:char} with itself and we obtain    
  \begin{eqnarray*}
    \c^{(3,1)\otimes(3,1)}=(1,0,1,1,9)
  \end{eqnarray*} 
  Since Table \ref{tab:char}  contains the row-vectors of all
  possible irreducible characters and that $(1,0,1,1,9)$ is
  obviously not one of these rows, it is immediate that the character
  represented by $(1,0,1,1,9)$ is not irreducible. But we observe
  that  the identity 
  $ 
    \c^{(3,1)\otimes(3,1)} = 
    \c^{(4)}+\c^{(3,1)}+\c^{(2,1,1)}+\c^{(2,2)}  
  $ 
  is true by adding the rows of Table \ref{tab:char} corresponding to
  the partitions in the right hand side.
  \hfill $\diamond$
\end{example}

More generally, the problem of finding the coefficients
$t_{\l,\m}^{\a}$ has a solution when one accepts to use the
character table $[\c^\l_\m]$ of $\S_n$:  
\begin{eqnarray} 
  \label{eq3}
  t_{\l,\m}^{\a}=\c^\l\otimes \c^\m|_{\c^\a} = 
  \sum_{\g\vdash n}\frac{|C_\g|}{n!}\c^\l_\g\c^\m_\g\c^\a_\g  
\end{eqnarray} 
Identity (\ref{eq3}) follows from the orthonormality of the characters
$\c^\l$ with respect to the standard scalar product in the group
algebra of $\S_n$ and from (\ref{eq1}).   But since the coefficients
$t_{\l,\m}^{\a}$ are positive integers, we find equation (\ref{eq3})
unsatisfactory and the goal of this paper is to contribute to other
avenues for computing the coefficients $t_{\l,\m}^{\a}$.   

Our  presentation is in three steps.  First we define
operators on the ring of symmetric functions  which reproduce tensor
product of irreducible representations when they act on Schur
functions.  Second we specialize these operators to the tensor powers
${\c^{(n-1,1)}}^{\otimes k}$, for $k\geq 0$, and we develop a
combinatorial model to represent these  tensor powers. Finally we
obtain some enumerative results for the multiplicity of any
irreducible representation in the  Kronecker powers of the form
$\c^{(n-1,1)\otimes k}$.    

\section{Combinatorial operators}
\label{sec:operators}

\subsection{Symmetric functions} 
\label{sub:symfunc}

Let $\mathbb{Q}[\S_n]$ be the group algebra of the symmetric group
over the field $ \mathbb{Q}$  of rational numbers and let
$\mathcal{}{Z}_n$ be the center of this group algebra. The irreducible
characters of $S_n$ can be seen as elements of $\mathcal{Z}_n$ when
one writes
\begin{eqnarray*}
  \c^\l=\sum_{\s\in\S_n}\c^\l(\s)\s
\end{eqnarray*} 
and the pointwise multiplication of two elements
$a=\sum_{\s\in\S_n}a_\s\s$ and $b=\sum_{\s\in\S_n}b_\s\s$ of
$\mathbb{Q}[\S_n]$ is defined by 
\begin{equation*}
  a\cdot b=\sum_{\s\in \S_n}(a_\s b_\s)\s
\end{equation*}
In example \ref{ex1} we have seen that pointwise multiplication of
characters gives the character of Kronecker product of the
corresponding representations: $\c^\l\cdot\c^\m=\c^{\l\otimes\mu}$. 

Let ${\x}=\{x_1,x_2,\ldots\}$ be a set of indeterminates,
$\L=\L_\mathbb{Q}[\x]$ the ring of symmetric functions in
$x_1,x_2,\ldots$ over the field $\mathbb{Q}$  and $\L^n$ the
restriction to homogeneous symmetric  functions of degree $n$. 
Two important sets of symmetric functions are the homogeneous
symmetric functions and the Schur symmetric functions. 
Given a partition $\l = (\l_1, \ldots, \l_m)$, one defines
$h_\l(\x)=\prod_{i=1}^m h_{\l_i}(\x)$, where $h_r(\x)$ is the $r$'th
homogeneous symmetric function, and one denotes by $s_\l(\x)$ the
Schur symmetric function associated to $\l$ (see \cite{macdo}). 
The sets $\{h_\l(\x)\}_{\l\vdash n}$ and $\{ s_\l(\x)\}_{\l\vdash n}$
are linear basis of $\L^n$. The Frobenius map $\mathcal{F}:
\mathcal{Z}_n\rightarrow \L^n$ is a vector space isomorphism which
sends the irreducible characters $\c^{\l}$ to the Schur functions
$s_\l$ : $\mathcal{F}(\c^{\l})=s_\l$.  
Schur functions can be expanded as determinants of homogeneous
functions: 
\begin{equation}
  \label{eq4}
  s_\l=\det(h_{\l_i-i+j})_{1\leq i,j\leq n}
\end{equation}
where $h_0=1$ and $h_r=0$ if $r<0$. The Littlewood-Richardson
coefficients (denoted by $LR$) are defined as the coefficients in the 
expansion of the ordinary product of two or more Schur functions in
the basis of Schur functions:  
\begin{equation*}
  s_{\nu^{(1)}}s_{\nu^{(2)}}\cdots s_{\nu^{(r)}} = 
  \sum_\m LR_{\nu^{(1)},\nu^{(2)},\ldots,\nu^{(r)}}^\m s_\m 
\end{equation*}
The adjoint operator to multiplication by $s_\g$ in $\L$ with respect
to the standard scalar product in $\L$  is denoted 
$s_\g^\perp$ and its action on the Schur function $s_\l$ is
 
as follows:
\begin{eqnarray}
  \label{eq5}
  s_\g^\perp s_\l &=& s_{\l/\g}= 
  \left\{ 
    \begin{array}{ll}
      \sum_\a LR_{\g,\a}^{\l}s_\a & \mbox{if $\gamma \subseteq \lambda$}\\
      0 & \mbox{otherwise}
    \end{array}
  \right.\\
  \label{eq6}
  \langle s_\g f,g \rangle &=&
  \langle f,s_\g^\perp g \rangle \qquad \forall \, f,g\in \L
\end{eqnarray} 

Now let us define  in $\L^n$ the operation $f\odot g$ corresponding to
pointwise multiplication by means of the Frobenius map:   
\begin{equation*}
  f\odot g = 
  \mathcal{F} (\mathcal{F}^{-1} (f)\cdot \mathcal{F}^{-1} (g)) \qquad
  \forall \,f,\,g,\, \in \L^n 
\end{equation*}

We shall call this operation\emph{ inner product} of symmetric
functions and we have in particular 
\begin{equation}
  \label{eq7}
  s_\l\odot s_\m=\sum_\a t_{\l,\m}^\a s_\a
\end{equation}
where the coefficients $t_{\l,\m}^\a $ are the same than in equation
(\ref {eq3}).  We can expand the inner product $h_\lambda\odot s_\mu$
in the basis \{$s_\a$\} as follows (see also \cite[1.7 example 23
(d)]{macdo}):  
\begin{eqnarray}\label{eq8}
  h_\lambda\odot s_\mu ({\mathbf x}) & 
  =                                  & 
  \sum_{\{ \nu^{(1)}\vdash \lambda_1,\nu^{(2)}\vdash\lambda_2,\ldots,
    \nu^{(k)}\vdash \lambda_k\}} LR_{\nu^{(1)}\nu^{(2)}\cdots
    \nu^{(k)}}^\mu 
  s_{\nu^{(1)}} \cdots  s_{\nu^{(k)}} \\ 
  \label{eq9}
    &
  = &
  \sum_{\{\nu^{(2)}\vdash\lambda_2,\ldots, \nu^{(k)}\vdash
  \lambda_k\}}(s_{\nu^{(2)}}\cdots
  s_{\nu^{(k)}})(s_{\nu^{(2)}}^\perp\cdots s_{\nu^{(k)}}^\perp )
  s_\mu. 
\end{eqnarray}

To prove  identity (\ref{eq8})  in the language of $\l$-rings, it
suffices to notice that $(h_\lambda\odot s_\mu) [{\mathbf X}] =
h_\l[{\mathbf {X  Y}}]|_{s_\m[{\mathbf Y}]}$. Then (\ref{eq8}) follows
from the fact that
$$
h_\l[\mathbf {XY}] = 
\sum_\m\left(\sum_{\{\nu^{(1)}\vdash\l_1,\ldots,\nu^{(k)}\vdash\l_k \}}  
  LR_{\nu^{(1)}\cdots \nu^{(k)}}^\mu s_{\nu^{(1)}} [{\mathbf X}]
  \cdots  s_{\nu^{(k)}}[{\mathbf X}]\right) s_\m[{\mathbf Y}]
$$ 
Identity (\ref{eq9}) follows from (\ref{eq8}) and  (\ref{eq5}).

\subsection{The operators $U_{\overline{\l}}$}
\label{sub:op-Ul}

\begin{definition}\em 
  Let $\l=(\l_1, \ldots, \l_m)\vdash n$ be a partition of $n$, with
  $\l_1\geq \l_2\geq \ldots \geq\l_m$, and $\overline{\l}$ the
  truncated partition of  $\l$ defined by
  $\overline{\l} = (\l_2, \ldots, \l_m)$.
  One denotes by  $U_{\overline{\l}}$ the operator from $\L^n$ to
  $\L^n$ defined as follows:
  \begin{description}
  \item[a)] 
    Expand the determinant 
    \begin{equation} \label{det}
      \left|\begin{array}{cccc}
          1 & 1 & \ldots & 1 
          \\
          h_{\l_2-1} & h_{\l_2} & \ldots & h_{\l_2+r-2} 
          \\
          \vdots & \vdots & \ldots & \vdots 
          \\
          h_{\l_r-r+1} & h_{_r-r+2}\ldots & \ldots & h_{\l_r}
        \end{array}\right|
    \end{equation}
  \item [b)]
    Replace each term $h_\a=h_{\a_1}h_{\a_2}\cdots h_{\a_m}$ in the
    expansion of  (\ref{det})  by  
    \begin{equation*}
      \sum_{\nu^1\vdash \a_1,\ldots , \nu^m\vdash \a_m}
      (s_{\nu^1}\cdot s_{\nu^m}) (s_{\nu^1}^\perp\cdots
      s_{\nu^m}^\perp) 
    \end{equation*}
    to obtain the operator $U_{\overline{\l}}$.
  \end{description}
\end{definition}

\begin{theorem}
  For any partitions $\l$  and $\mu$ of $n$,
  \begin{eqnarray*}
    U_{\overline{\l}} s_\m &=& s_\l\odot s_\m 
    \\
    &=& \sum_{\a\vdash n}t_{\l,\m}^\a s_\a
  \end{eqnarray*}
\end{theorem}

\begin{proof}
  This is a straightforward consequence of equations  (\ref{eq4} ),
  (\ref{eq7}) and (\ref{eq9}).
\end{proof}

\begin{example}\label{ex2.1}\em
  The Kronecker product $\c^{(n-1,1)}\otimes 
  \c^\m$ is obtained by applying the operator $U_{(1)}$ on $s_\m$
  which is obtained by  expanding the determinant in definition 2.1 a)
  and then writing the $h_\l$ in terms of Schur functions using
  definition 2.1 b):  
  \begin{equation*}
    \left|\begin{array}{cc}
        1 & 1 
        \\
        h_0 & h_1
      \end{array}\right|
    =
    h_1-1\Rightarrow U_{(1)}=s_{(1)}s_{(1)}^\perp-1
  \end{equation*}
  Similarly the operator $U_{(2)}$ needed for the computation of
  $\c^{(n-2,2)}\otimes  \c^\m$ is obtained with the determinant
  \begin{equation*}
    \left|\begin{array}{cc}
        1 & 1 
        \\
        h_1 & h_2
      \end{array}\right|
    =
    h_2-h_1\Rightarrow U_{(2)}=s_{(2)}s_{(2)}^\perp +
    s_{(1,1)}s_{(1,1)}^\perp-s_{(1)}s_{(1)}^\perp
  \end{equation*}
  ~\hfill $\diamond$
\end{example}
                
The fact that the computation of $\c^{\l}\otimes\c^\m$ is independent
of the  largest part of $\l$  was already observed by Murnaghan (\cite{Mu2}) and
also described by Thibon in \cite{Th2}, but the definition of the
operators  $U_{\overline{\l}}$ seems to be new.  

\section{
  A combinatorial model for $\c^{(n-1,1)\otimes k}$ and some
  consequences
}  

We are now ready to present a combinatorial model for the multiplicity 
of any irreducible representation in any Kronecker power
$\c^{(n-1,1)\otimes k}$, in terms of paths in Young's lattice.  
  From example \ref{ex2.1}, it is immediate that  the expansion of the 
  Kronecker power $\c^{(n-1,1)\otimes k}$ is obtained by the
  application of the operator ${U_{(1)}}^k=(s_{(1)}s_{(1)}^\perp-1)^k$ on the
  Schur function $s_{(n)}$.
In the remaining of this section we develop a combinatorial
interpretation of this process, followed by some enumerative
consequences. 

Let $\l = (\l_1, \ldots, \l_m)$ be an integer partition of $n$.
We recall that the unique {\em Ferrers diagram} associated to $\l$ is
formed of $m$ stacked rows of cells, ordered from bottom to top in
increasing order (in the french notation) and such that the $i^{th}$
row contains $\l_i$ cells.
A cell of a diagram $\l$ located at the right end of a row and having no
cell above it is called a {\em corner} of $\l$.  
For example, in the following Ferrers diagram, corresponding to the
partition $(4,4,2,1)$, the corners are indicated by $\bullet$. 
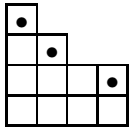
$$
  \young(\bullet,\nil\bullet,\nil\nil\nil\bullet,\nil\nil\nil\nil)
$$ 
The unique corner of a Ferrers diagram located on a longest row is
called the {\em first corner} of the diagram.   

It follows immediately from the definition of the operator
$s_{(1)}^\perp$ that, for a Ferrers diagram $\l$, $s_{(1)}^\perp(s_\l)$ is 
the sum of the Schur functions indexed by the Ferrers diagrams
obtained by removing a single corner from $\l$.  
Symmetrically, $s_{(1)}(s_\l)$ is the sum of the Schur functions indexed
by the Ferrers diagrams obtained by adding to $\l$ a new cell that
becomes a corner of the new diagram. 
Hence $s_{(1)}s_{(1)}^\perp(s_\l)$ is the sum of the Schur functions
indexed with the  Ferrers diagrams that are obtained from $\l$ by
first removing a corner from $\l$, which gives a diagram denoted
$\l'$, then adding a corner to $\l'$.  
One says that every diagram, or equivalently partition, $\l' \neq \l$
indexing a Schur function occurring in the sum $s_{(1)}
s_{(1)}^\perp(s_\l)$ {\em differs from $\l$ by the position of a
  corner}.        
It should be noticed that in the sum $s_{(1)}s_{(1)}^\perp(s_\l)$ the
multiplicity of $s_\l$ is at least $1$.    
Hence, as $U_{(1)}=s_{(1)}s_{(1)}^\perp - 1$, one can define $U_{(1)}$
as $s_{(1)}s_{(1)}^\perp(s_\l)$ minus one occurrence of $s_\l$.    

\begin{example}\label{ex3.1}\em
  $U_{(1)}(s_{(3,3,1)}) = s_{(3,3,1)} + s_{(4,2,1)} + s_{(3,2,1,1)} +
  s_{(3,2,2)} + s_{(4,3)}$, and the partitions $(4,2,1)$, $(3,2,1,1)$,
  $(3,2,2,2)$ and $(4,3)$ differ from $(3,3,1)$ by the position of a
  corner. The 
  multiplicity of $s_{(3,3,1)}$ is $1$ because there are two corners
  in $(3,3,1)$, and thus only two ways to obtain $(3,3,1)$ from itself
  by removing then replacing a corner, one of these occurrences being
  not tamen into account due to the $-1$ term in the definition of
  $U_{(1)}$. 
  \hfill
  $\diamond$
\end{example}

We now define the main combinatorial object that we will need in order
to encode the action of $U_{(1)}$ on a Schur function $s_\m$ iterated
$k$ times. 

\begin{definition}\em 
  \label{def:OST}
  Given a positive integer $k$ and partitions $\l$ and $\m$ of same weight, a {\em
    Kronecker tableau} $K$ of length $k$, initial shape $\mu$ and
  final shape $\l$ is a sequence $\mu^0=\mu, \mu^1, \ldots,
  \mu^k=\l$ of Ferrers diagrams where, 
 for every pair of consecutive diagrams $\mu^i$ and $\mu^{i+1}$,
  {\it 
    either $\mu^{i+1}$ differs from $\mu^i$ by the position of a
    corner,  
    or $\mu^{i+1}=\mu^i$ and one corner of $\mu^{i+1}$, other than its
    first corner, is distinguished. 
  } 
  One denotes by $KT_{\mu,\l}^k$ the set of {\em Kronecker tableaux} 
  of length $k$, initial shape $\m$ and final shape $\l$. 
\end{definition}

\begin{example}\em
  \label{ex:OST}
  The following Kronecker tableau -- where a distinguished corner is 
  indicated by a $\times$ -- belongs to $KT^{9}_{(5),(3,2)}$
  $$
  \Yboxdim8pt
  \yng(5) \yng(1,4) \young(\times,\nil\nil\nil\nil) \yng(2,3)
  \yng(1,1,3)  \yng (1,2,2)
  \young(\times,\nil\nil,\nil\nil) \yng(1,1,1,2) \yng(1,2,2)
  \yng(2,3)  
  \Yautoscale1
  $$
\end{example}

\begin{proposition}
  \label{prop:OST-model}
  Let $k$ and $n$ be two positive integers and $\l$ a partition of
  $n$. Then  
  \begin{equation}
    \label{eq:OST-model}
    {\c^{(n-1,1)}}^{\otimes k}|_{\c^\l}
    = 
    |KT_{{(n)},{\l}}^k|.
  \end{equation}
\end{proposition}
\begin{proof}
  The definition of $U_{(1)}$ as an operator on Schur functions can be 
  translated in the combinatorial framework of partitions and Ferrers
  diagrams, due to the fact that Schur functions are indexed by
  partitions. 
  Hence, the sum of Schur functions ${U_{(1)}}^k(s_n)$ can be seen as a
  formal sum of Ferrers diagrams.
  The number of occurrences of a diagram $\l$ in this formal sum of
  diagrams is then given by the number of ways to obtain $\l$ from
  $(n)$ by iterating $k$ times the combinatorial operation associated
  to $U_{(1)}$.
  The identity follows immediately from this fact and from the
  definition of Kronecker tableaux, where the restriction that the
  first corner of a diagram can not be distinguished in the next
  diagram accounts for the $-1$ part of $U_{(1)}=s_{(1)}s_{(1)}^\perp-1$.
\end{proof}

The above proposition establishes a link between the multiplicity of the
irreducible character $ \chi^\l$ in some Kronecker  power
and sequences of Ferrers diagrams  seen as paths in Young's lattice
(the lattice of Ferrers diagrams ordered by inclusion). 
We rely on this fact to obtain enumerative results about the
multiplicity of any given irreducible representation in the Kronecker  
power $\c^{(n-1,1)\otimes k}$.
The main tool we use is a combinatorial construction defined for a
family of paths in Young's lattice called {\em oscillating tableaux}
and introduced by Sundaram, in a different algebraic context, in
\cite{Su} (see also the work of Delest, Dulucq and Favreau
\cite{DDF,Fa} for a purely combinatorial point of view). 

Briefly, oscillating tableaux are paths in Young's lattice, that is
sequences of Ferrers diagrams, starting at $\emptyset$ and such that
{\em two consecutive diagrams differ by the addition or removal of
  exactly one corner}. 
\begin{example}\em
  \label{ex:OT}
  Here is an oscillating tableau of length $7$ and final shape
  $(2,1)$, that contains five additions of corner and two removal of 
  corner (steps $4$ and $7$).
  $$
   \Yboxdim8pt
  \emptyset \yng(1) \yng(1,1) \yng(1,2) \yng(2) \yng(3) \yng(1,3)
  \yng(1,2) 
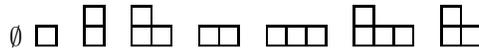
  $$
  \hfill $\diamond$
\end{example}

In Lemmas \ref{lem:bijection1} and \ref{lem:bijection2} below, we
consider a class of Kronecker tableaux that can be related to
oscillating tableaux. This allows to use a variant of the
combinatorial construction defined in \cite{Su,DDF} that will be
central in the proof of our main enumerative result, Proposition
\ref{prop:enumeration}.

\begin{lemma}
  \label{lem:bijection1}
  Let $n$ and $k$ be positive integers and $\l$ a partition of 
  $n$ such that $n \geq k + \l_2$. 
  There is a bijection between Kronecker tableaux of $KT_{(n),\l}^k$
  and sequences $\n^0, \ldots, \n^k$ of $k$ Ferrers diagrams such that
  $\n^0=\emptyset$, $\n^k=\overline{\l}$ and, for every pair $\n^i$
  and $\n^{i+1}$ of consecutive diagrams, either $\n^{i+1}$ is obtained
  from $\n^i$ by the addition or removal of one corner, or $\n^{i+1}$
  differs from $n^i$ by the position of a corner, or $\n^{i+1}=\n^i$
  and $\n^{i+1}$ has one distinguished corner. 
\end{lemma}
\begin{proof}
  Let $K$ be a Kronecker tableau of $KT_{(n),\l}^k$ such that $n \geq
  k + \l_2$.  
  Due to this last condition, the first corner of every Ferrers
  diagram of $K$, except possibly the last diagram,  is on its first row.
  Then, by removing the first row of every diagram of $K$ one obtains
  the sequence $\n^0, \ldots, \n^k$. 
  Conversely, consider a sequence of $k+1$ Ferrers diagrams
  $\n^0=\emptyset, \ldots, \n^k=\overline{\l}$. 
  By adding, for every diagram $\n^i$, a first row of length
  $n-|\n^i|$, one obtains a Kronecker tableau of $KT_{(n),\l}^k$. 
\end{proof}

The combinatorial construction we describe below relies partly on the
{\em Robinson-Schensted-Knuth (RSK) insertion and deletion
  algorithms}, and we first recall some basic facts about {\em
  standard tableaux} and these algorithms (see \cite{Sa} for
example for details on these algorithms).
\begin{itemize}
\item
  Given a positive integer $n$ and a Ferrers diagram $\alpha$ with at 
  most $n$ cells, a {\em partial standard tableau} of shape $\alpha$
  and labels in $[n]$ is a labelling of the cells of $\alpha$ with
  distinct integers chosen from $\{1, \ldots, n\}$, such that the
  labels are increasing in rows (from left to right) and columns
  (from bottom to top).  
\item
  Let $S$ be a partial standard tableau. Given an integer $x$, the
  RSK insertion algorithms inserts $x$ into $S$, creating a tableau
  $S'$ whose shape differs from the shape of $S$ by the addition of a
  corner and labels are the labels of $S$ plus $x$. Given a corner of
  $S$, the RSK deletion algorithm removes this corner and moves some
  labels of cells of $S$, this process ending when a label is ejected
  from the first row of the tableau. 
\end{itemize}

\begin{lemma}
  \label{lem:bijection2}
  Let $n$ and $k$ be two positive integers and $\l$ a partition of
  $n$ such that $n \geq k + \l_2$.
  There is a bijection between the set $KT^k_{(n),\l}$ and the
  set of pairs  $(T,\pi)$, where
  $T$ is a partial standard tableau of shape $\overline{\l}$ with
  labels in $[n]$ and $\pi$ is a permutation of $n$ such that: every
  non fixed point cycle of $\pi$ is decreasing, every fixed point of
  $\pi$ is also the label of a cell of $T$, and every label of $T$ is
  the greatest element of a cycle in $\pi$.
\end{lemma}

\begin{proof}
  Let $K=\m^0, \ldots, \m^k$ be a Kronecker tableau of length $k$,
  initial shape $(n)$ and final shape $\l$, such that $n \geq k +
  \l_2$. Let $\n^0, \ldots, \n^k$ be the sequence of Ferrers diagrams
  corresponding to $K$ obtained by the construction of Lemma
  \ref{lem:bijection1}.

  One can associate to $\n^0, \ldots, \n^k$ a sequence of partial
  standard tableaux $(T_0=\emptyset, \ldots, T_k=T)$ on $[n]$ and a
  permutation $\pi$, such  that the shape of $T_i$ is $\n^i$ for every
  $i$. We proceed  as follows. Start with setting $\pi$ as the
  identity permutation on $\{1,\ldots ,k\}$, and for $i$ from $1$ to
  $k$: 
  \\ {\em 1.}  
  If $\n^{i}$ is obtained from $\n^{i-1}$ by the addition of a corner,
  then add to $T_{i-1}$ this corner, labelled with $i$, to obtain
  $T_i$.
  \\ {\em 2.}
  If $\n^{i}$ is obtained from $\n^{i-1}$ by the removal of a corner,
  then delete this corner from $T_{i-1}$, using the RSK deletion
  algorithm. If $j$ is the integer ejected from $T_{i-1}$ by the
  RSK deletion algorithm, then multiply $\pi$ by the transposition
  $(i,j)$. 
  \\ {\em 3.} 
  If $\nu^i$ differs from $\nu^{i-1}$  by the position of a corner,
  or $\nu^i=\nu^{i-1}$ and $\nu^{i-1}$ has a distinguished corner
  (therefore $\nu^i$ and $\nu^{i-1}$ have the same weight), then
  delete this corner from $T_{i-1}$ using again the RSK deletion
  algorithm and add the corner needed to obtain $\nu^i$ and label it
  with $i$.   
  If $j$ is the ejected label  then multiply $\pi$ by the
  transposition $(i,j)$.  

  The fact that in the permutation $\pi$ all non fixed point cycles 
  are decreasing follows from the fact that in every transposition
  $(i,j)$ considered in the construction above one has $i > j$.

  The reverse construction starts with a  partial standard tableau
  $(T_k=T)$  with shape $\n^k$ and a permutation $\pi$ of the set
  $\{1,\ldots , k\}$  with each non fixed point cycle in decreasing
  order such that each entry of $T$ is the greatest element of a cycle 
  of $\pi$ (including the fixed points).   
  Then perform the following steps, for $i$ from $k$ to $1$:
  \\ {\em 1.}
  If no cell of $T_i$ is labelled with $i$, there exists $j < i$ such
  that $\pi(i)=j$. Then insert the integer $j$ into the tableau $T_i$
  using the RSK insertion algorithm to obtain $T_{i-1}$ and define 
  $\n^{i-1}$ as the shape of $T_{i-1}$. 
  \\ {\em 2.}
  If a cell of $T_i$ is labelled with $i$, then remove the cell
  labelled $i$: by induction this cell is a corner and this removal 
  gives a partial standard tableau denoted $U$.
  \\ {\em 2.a.}
  If furthermore there exists $j < i$ such that $\pi(i)=j$, then
  insert the integer $j$ into the tableau $U$, using the RSK insertion 
  algorithm, to obtain $T_{i-1}$, and define $\n^{i-1}$ as the shape
  of $T_{i-1}$, distinguishing the corner added during this insertion 
  if it takes the same position than the corner removed from $T_i$.
  \\ {\em 2.b.}
  Otherwise, after removing $i$ from $T_i$, multiply $\pi$ by the
  transposition $(i,j)$, and define $\n^{i-1}$  as the shape of
  $T_{i-1}$.  
  
  The fact that these two constructions define  a bijection follows
  immediately from its close relationship with the construction on
  oscillating tableaux defined by Sundaram \cite{Su} and Delest,
  Dulucq and Favreau \cite{DDF,Fa}.  
\end{proof}

\begin{example}\em
  \label{ex:bijection}
  The following Kronecker tableau belonging to $KT^{12}_{(6),(2,2,2)}$
  $$
  \Yboxdim8pt
  \yng(6) \yng(1,5) \young(\times,\nil\nil\nil\nil\nil) \yng(2,4)
  \yng(1,2,3) \yng(1,1,4)  \yng (1,2,3) \yng(2,2,2) \yng(1,1,2,2)
  \yng(1,2,3) \yng(2,2,2) \yng(1,2,3) \yng(2,2,2)
  \Yautoscale1
  $$
  corresponds to the sequence of partial standard tableaux $\n^0,
  \ldots, \n^k$   
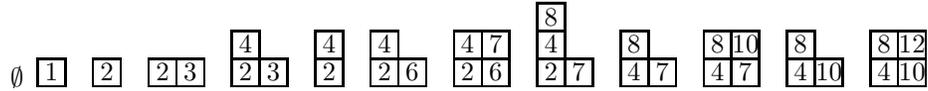
  $$
  \Yboxdim11pt
  \emptyset \young(1) \young(2) \young(23) \young(4,23) \young(4,2)
  \young(4,26)  \young(47,26) \young(8,4,27) \young(8,47)
  \young(8\ten,47)  \young(8,4\ten) \young(8\twelve,4\ten) 
  \Yautoscale1
  $$
  and to the pair
  $$
  \Yvcentermath1
  T = \young(8\twelve,4\ten),\  
  \pi =(1,2,\ldots12).(2,1).(5,3).(8,6).(9,2).(11,7)
  =(4)(5,3)(8,6)(9,2,1)(10)(11,7)(12).
  $$
  \hfill $\diamond$
\end{example}

To conclude this section, we derive from the above bijection an
explicit formula and a generating function for the coefficients
$\c^{(n-1,1)})^{\otimes k}|_{\c^\l}$ when $n\geq k+\l_2$.

\begin{proposition}
  \label{prop:enumeration}
  Let $k$ and $n$ be two positive integers and $\l$ a partition of
  $n$ such that $n\geq k+\l_2$.  Then
  \begin{equation}
    \label{eq:enumeration}    
    {\c^{(n-1,1)}}^{\otimes k}|_{\c^\l} =
    f^{\overline{\l}}
    \sum_{m_1=0}^{n-\l_1}
    \left(
      \binom{k}{m_1}
      \sum_{m_2=n-\l_1-m_1}^{\lfloor (k-m_1)/2 \rfloor}
    \binom{m_2}{n-\l_1-m_1} p_2(k-m_1,m_2) 
    \right),
  \end{equation}
  where $f^{\overline{\l}}$ is the number of standard tableaux of shape
  $\overline{\l}$ and $p_2(k-m_1,m_2)$ is the number of set partitions
  of a set of $k-m_1$ distinct integers into $m_2$ parts of size at
  least $2$.
\end{proposition}

\begin{proof}
  From Proposition \ref{prop:OST-model}, it is sufficient to enumerate 
  the number of Kronecker tableaux of length $k$, initial shape $(n)$  and final shape
  $\l$. As $n\geq k+\l_2$, it follows from Lemma \ref{lem:bijection2}
  that this reduces to the enumeration of some couples $(T,\pi)$ where
  $T$ is a partial standard tableau of shape $\overline{\l}$ and $\pi$ is a
  permutation on the set $\{1,\ldots ,k\}$. 
 Formula (\ref{eq:enumeration}) follows if one denotes by $m_1$ the
  number of fixed points of $\pi$ and $m_2$ the number  of cycles of
  size at least $2$ in $\pi$. 
\end{proof}

\begin{remark}\em
  \label{rem:p2} 
  For two given integers $n$ and $k$, the number $p_2(n,k)$ is known
  as an associated Stirling number of second kind (reference A008299
  in \cite{Sloane}, see also \cite[p. 222]{Comtet}). 
  Such numbers are defined by the following recurrence: $p_2(n,k) = 0$
  if $n < 2k$ and $p_2(n,k) = kp_2(n-1,k)+ (n-1)p_2(n-2,k-1)$ if $n \geq
  2k$. The computation of $p_2(n,k)$ can also be done by extracting
  the coefficient of $y^k x^n/n!$ in $e^{y p(x)}$ where $p(x) =
  e^x-x-1$. 
\end{remark}

\begin{corollary}
  \label{cor:gf}
  Let $\ell$ be a positive integer, $\overline{\l}=(\ell_2,\ldots ,
  \ell_m)$ an integer partition of $\ell$. Then  
  \begin{equation}
    \label{eq:gf-main}
    \sum_{k \geq \ell} 
    {\c^{(n_k-1,1)}}^{\otimes k}| _{\c^{{\overline{\l}}+(n_k-\ell)}}
    \frac{x^k}{k!} 
    =  
    \frac{f^{{\overline{\l}}}} {\ell!}  e^{p(x)} (e^x-1)^{\ell},
  \end{equation}
  where, for every $k \geq \ell$,  $n_k$ is any integer satisfying
  $n_k\geq k+\l_2$ and
  ${\overline{\l}}+(n_k-\ell)$ is the integer partition  obtained by
  adding the part $n_k-\ell$ to ${\overline{\l}}$ . 
\end{corollary}

\begin{proof}
  For clarity, we denote
  $\overline{\l}^*=({\overline{\l}}+(n_k-\ell))$  
  and  ${p(x)}=e^x-x-1$. 
  It follows from the fact that $n_k\geq k+\lambda_2$, lemma
  \ref{lem:bijection1}, propositions \ref{prop:OST-model} and
  \ref{prop:enumeration}, that 
  \begin{eqnarray*}
   \sum_{k \geq \ell} 
   {\c^{(n_k-1,1)}}^{\otimes k}|_{\c^{(\overline{\l}^*)}} \frac{x^k}{k!} 
   & = &
    f^{{\overline{\l}}}
    \sum_{k \geq \ell} 
    \sum_{m_1=0}^\ell 
    \left(
      \frac{x^k}{k!}
      \binom{k}{m_1}
      \sum_{m_2=\ell-m_1}^{\lfloor (k-m_1)/2\rfloor } 
      \binom{m_2}{\ell-m_1}
      p_2(k-m_1,m_2) 
          \right)
    \\
    \Longleftrightarrow
    \sum_{k \geq \ell} {\c^{(n_k-1,1)}}^{\otimes
    k}|_{\c^{(\overline{\l}^*)}} \frac{x^k}{k!}  
    & = &
    f^{{\overline{\l}}}
    \sum_{m_1 = 0}^\ell
    \frac{x^{m_1}}{m_1!}
    \left(
      \sum_{k \geq \ell}
      \sum_{m_2=\ell-m_1}^{\lfloor(k-m_1)/2\rfloor} 
      \binom{m_2}{\ell-m_1}
      p_2(k-m_1,m_2) 
      \frac{x^{k-m_1}}{(k-m_1)!}
    \right)
    \\
    \Longleftrightarrow
    \sum_{k \geq \ell} {\c^{(n_k-1,1)}}^{\otimes
      k}|_{\c^{(\overline{\l}^*)}} \frac{x^k}{k!}  
    & = &
    f^{{\overline{\l}}}
    \sum_{m_1 = 0}^\ell
    \frac{x^{m_1}}{m_1!}
    \left(
      \sum_{m_2 \geq \ell-m_1}
      \binom{m_2}{\ell-m_1}
      \sum_{q \geq 0}
      p_2(q,m_2)
      \frac{x^{q}}{q!}
    \right)
    \\
    \Longleftrightarrow
    \sum_{k \geq \ell} {\c^{(n_k-1,1)}}^{\otimes
      k}|_{\c^{(\overline{\l}^*)}} \frac{x^k}{k!}  
    & = &
   f^{{\overline{\l}}}
    \sum_{m_1 = 0}^\ell
    \frac{x^{m_1}}{m_1!}
    \left(
      \sum_{m_2 \geq \ell-m_1}
      \binom{m_2}{\ell-m_1}
      \frac{p(x)^{m_2}}{m_2!}
    \right)
    \\
    \Longleftrightarrow
    \sum_{k \geq \ell} {\c^{(n_k-1,1)}}^{\otimes
      k}|_{\c^{(\overline{\l}^*)}} \frac{x^k}{k!}  
    & = &
   f^{{\overline{\l}}}
    \sum_{m_1 = 0}^\ell
    \frac{x^{m_1}}{m_1!}
    \frac{p(x)^{\ell-m_1}}{(\ell-m_1)!}
    \left(
      \sum_{m \geq 0}
      \frac{p(x)^{m}}{m!}
    \right)
      \\
    \Longleftrightarrow
    \sum_{k \geq \ell} {\c^{(n_k-1,1)}}^{\otimes
      k}|_{\c^{(\overline{\l}^*)}} \frac{x^k}{k!}  
    & = & 
   f^{{\overline{\l}}}
    \sum_{m_1 = 0}^\ell
    \frac{x^{m_1}}{m_1!}
    \frac{p(x)^{\ell-m_1}}{(\ell-m_1)!}
    e^{p(x)}
    \\
    \Longleftrightarrow
    \sum_{k \geq \ell} {\c^{(n_k-1,1)}}^{\otimes
      k}|_{\c^{(\overline{\l}^*)}} \frac{x^k}{k!}  
    & = &
    \frac{f^{{\overline{\l}}}}{\ell!} e^{p(x)} (e^x-1)^{\ell}.
  \end{eqnarray*}  
\end{proof}

\section{Conclusion}
We presented in this note a combinatorial interpretation of the
multiplicity of any irreducible representation in a Kronecker power
${\c^{(n-1,1)}}^{\otimes k}$ in terms of sequences of Ferrers diagrams 
(Kronecker tableaux) that leads, when $n$ is large enough with respect
to $k$ and $\overline{\l}$, to an enumeration formula and a generating
function. Moreover, we now have a combinatorial model for the
expansion of ${\c^\mu}^{\otimes k}$ for any $\mu$ given by the
differential operators $U_{\overline{\mu}}$. However, at this point,
the problem of transforming this model into enumerative results in
terms of a relationship between $n$, $k$, $\mu$ and $\l$ is still
open.

It should be noticed that our approach could be extended to the more
general case of the computation of ${\c^{(n-1,1)}}^{\otimes k} \otimes
\c^\mu |_{\c^\l}$ for arbitrary $\mu$. This would require the use of a
generalization of the oscillating tableaux of Sundaram which already
exists and is called {\em skew oscillating tableaux} in \cite{DuSa}
(see also \cite{ChDu}).  However, this construction leads to an
intricate expression for the enumeration of Kronecker tableaux of
initial shape $\mu$ with more than one part, and we were not able to
find a compact generating function similar to Corollary \ref{cor:gf}.


\end{document}